\documentclass[11pt,a4paper]{article}
\usepackage{latexsym,epsfig}
\usepackage{graphicx} 
\usepackage[usenames,dvipsnames]{color} 
\usepackage{url} 
\usepackage{soul} 
\usepackage{array} 
\usepackage{amssymb}
\usepackage{amsfonts}
\usepackage{amsmath}
\usepackage{verbatim}
\usepackage{hyperref}
\usepackage[ruled]{algorithm2e}

\newtheorem{theorem}{Theorem}[section]
\newtheorem{corollary}{Corollary}[section]
\newtheorem{problem}{Problem}

\newcommand{\EndProof}{\hspace{\stretch{1}} $\Box$}

\newcommand{\noi}{\noindent}

\newcommand{\NP}{\mathcal{NP}}

\title{The Maximum Degree\&Diameter-Bounded Subgraph in the Mesh}

\author{
Mirka Miller$^{1, 2, 3}$\footnote{E-mail: Mirka.Miller@newcastle.edu.au. This research was supported by a Marie Curie International Incoming Fellowship within the 7th European Community Framework Programme.} \\ \,
Hebert P\'erez-Ros\'es$^{1}$\footnote{E-mail: Hebert.Perez@gmail.com} \\ \,
Joe Ryan$^{1}$\footnote{E-mail: Joe.Ryan@newcastle.edu.au} \\ \\
$\phantom{0}^{1}$\emph{Department of Software Engineering and Computer Science}, \\
\emph{The University of Newcastle, Australia} \vspace{3mm} \\   
$\phantom{0}^{2}$\emph{Department of Informatics, King's College, London, UK} \\ 
$\phantom{0}^{3}$\emph{Department of Mathematics, University of West Bohemia, Pilsen, Czech Republic} }

\begin{document}
\maketitle

\thispagestyle{empty}
\begin{abstract}
\noi The problem of finding the largest connected subgraph of a given undirected host graph, subject to constraints on the maximum degree $\Delta$ and the diameter $D$, was introduced in \cite{maxddbs}, as a generalization of the Degree-Diameter Problem. A case of special interest is when the host graph is a common parallel architecture. Here we discuss the case when the host graph is a $k$-dimensional mesh. We provide some general bounds for the order of the largest subgraph in arbitrary dimension $k$, and for the particular cases of $k=3, \Delta = 4$ and $k=2, \Delta = 3$, we give constructions that result in sharper lower bounds.
\end{abstract}

\noi \textbf{Keywords:} Network design, Degree-Diameter Problem, parallel architectures, mesh, Delannoy numbers 
\vskip 5mm


\section{Introduction}
\label{sec:intro}

Let $G = (V, E)$ be an undirected graph without loops or multiple edges (called the \emph{host graph}), with $n$ vertices (its \emph{order}), and $m$ edges (its \emph{size}). Our problem is stated as follows:

\begin{problem}[\textsc{MaxDDBS}]
\label{prob:ddbs}
Given a connected undirected host graph $G$, an upper bound $\Delta$ for the maximum degree, and an upper bound $D$ for the diameter, find the largest connected subgraph $S$ with maximum degree $\leq \Delta$ and diameter $\leq D$.
\end{problem}

\textsc{MaxDDBS} is a natural generalization of the well-known Degree-Diameter Problem (DDP), which asks for the largest graph with given degree and diameter \cite{Moore-survey}. DDP can be seen as \textsc{MaxDDBS} when $G$ is the complete graph $K_n$ for sufficiently large $n$. Problem \ref{prob:ddbs} was recently introduced in \cite{maxddbs}, where various practical applications are discussed and a heuristic approximation algorithm to solve \textsc{MaxDDBS} is given, since it is computationally hard.

Regarding computational complexity, \textsc{MaxDDBS} is known to be $\NP-$hard, since it contains other well-known $\NP-$hard problems as subproblems. In fact, restricting the search to only one constraint (either on the degree or the diameter), is enough to ensure $\NP-$hardness \cite{Johnson85}. The Largest Degree-Bounded Subgraph Problem is $\NP-$hard as long as we insist that the subgraph be connected, but can be solved in polynomial time otherwise (Problem ND1 of \cite{Garey-Johnson}). On the other hand, the Maximum Diameter-Bounded Subgraph becomes the Maximum Clique for $D = 1$, which was one of Karp's original 21 $\NP-$hard problems \cite{Karp72}. \textsc{MaxDDBS} also turns out not to be in \textsc{Apx}, the class of $\NP-$hard optimization problems for which there is a polynomial-time algorithm with a constant approximation ratio.

As mentioned above, \textsc{MaxDDBS} is also closely related to the Degree-Diameter Problem (DDP), stated by Elspas in 1964, which consists of finding the largest graph with a given maximum degree $\Delta$ and a given diameter $D$. Since the order of such a graph cannot exceed the quantity $M_{\Delta,D} = 1 + \Delta + \Delta(\Delta -1) + \dots + \Delta(\Delta -1)^{D-1}$, called the \emph{Moore bound}, if we take $G$ as the complete graph on $M_{\Delta,D}$ vertices (denoted by $K_{M_{\Delta,D}}$) in Problem \ref{prob:ddbs}, we get the Degree-Diameter Problem. Note that this does not imply that DDP is $\NP-$hard; actually, the complexity of DDP is not known to-date. Obviously, the Moore bound is also a theoretical upper bound for \textsc{MaxDDBS}.

A graph whose order is equal to the Moore bound is called a \emph{Moore graph}. Moore graphs are very rare; they exist only for certain special cases; for diameter $D=1$, Moore graphs are the complete graphs of order $\Delta+1$, for maximum degree $\Delta = 2$, Moore graphs are the odd cycles. The only other Moore graphs are of diameter $D = 2$, $\Delta= 3, 7$ and possibly $57$, \cite{Moore-survey}. We denote by $N_{\Delta, D}$ the order of the largest graph that can be constructed with maximum degree $\Delta$ and diameter $D$; the current lower bounds for $N_{\Delta, D}$ are shown in \cite{combinatoricswiki}.

A case of special interest is when the host graph $G$ is a common parallel architecture, such as the mesh, the hypercube, the butterfly, or the cube-connected cycles. If there are any constraints on communication time between two arbitrary processors, then \textsc{MaxDDBS} corresponds to the largest subnetwork that can be allocated to perform the computation. The case of the mesh and the hypercube as host graphs were already treated in \cite{maxddbs}, where some bounds were found for the order of \textsc{MaxDDBS} in a $k$-dimensional mesh.
Here we revisit in more detail the case of the mesh as a host graph. We refine the bounds given in \cite{maxddbs} for the order of the largest subgraph in arbitrary $k \geq 1$, and we focus on the cases $k=3, \Delta=4$ and $k=2, \Delta=3$. For those particular cases we give constructions that result in larger lower bounds.

The rest of the paper is organized as follows: In Section \ref{sec:arbitrary} we give the bounds for \textsc{MaxDDBS} in a $k$-dimensional mesh. Then, in Section \ref{sec:3D} we give the constructions for $\Delta = 4$ in the $3$-dimensional mesh. Section \ref{sec:2D} gives constructions for $\Delta = 3$ in dimension two. Finally, in Section \ref{sec:open} we discuss some open problems and research directions.


\section{MaxDDBS in the $k$-dimensional mesh}
\label{sec:arbitrary}

\begin{figure}[htbp]
\begin{center}
	 	\includegraphics[width=0.6\textwidth]{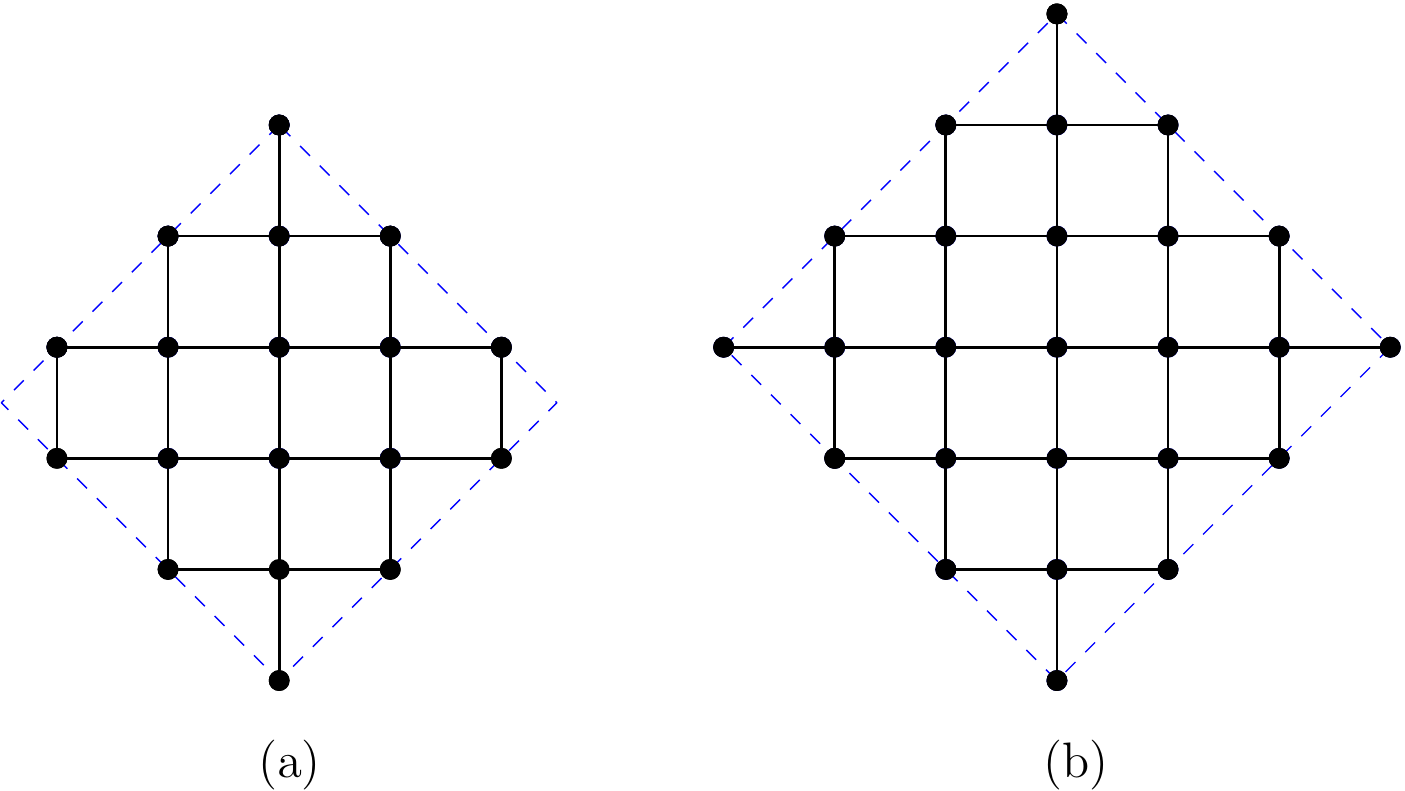}
	 \caption{Maximal balls in the two-dimensional $L^1$ metric space}
	 \label{fig:balls2D}
\end{center}	
\end{figure}

Here we will assume that the host graph $G$ is an infinite $k$-dimensional mesh, and we are looking for a subgraph of maximum degree $\Delta \leq 2k$, and diameter $D$. We can associate our mesh with an $L^1$ metric space in dimension $k$. Pick an arbitrary point in this $k$-dimensional $L^1$ metric space as the center of a coordinate system. Now, the vertices of the mesh are the points with integer coordinates (\emph{lattice points}), and two lattice points are joined iff they are at distance $1$. The largest subgraph of degree $\Delta = 2k$ and diameter $D$ corresponds to a closed ball of radius $D/2$. The number of lattice points contained in a ball of radius $D/2$ is variable, and depends on the location of the center of the ball. The maximum number of lattice points is achieved when the center of the ball is a lattice point itself, for even $D$, and when the center of the ball is the midpoint between two adjacent lattice points, for odd $D$. Such balls will be called  \emph{maximal}, and will be denoted $B_k(p)$, where $p = \lfloor \frac{D}{2} \rfloor$, depending on its parity, and $k \geq 1$. Figure \ref{fig:balls2D} depicts two maximal balls in dimension two, with diameters $5$ and $6$, respectively.

For more details about the shape and location of maximal balls, we refer to the comprehensive study by Dougherty and Faber \cite{dougherty}, where the same problem appears in a slightly different context. In that paper, the infinite $k$-dimensional mesh is interpreted as the Cayley graph of the free abelian group $\mathbb{Z}^k$, and the points with integer coordinates contained in the ball $B_k(p)$ correspond to words of length at most $D$ in the canonical generators of $\mathbb{Z}^k$. The aim of \cite{dougherty} was to construct large Cayley graphs on abelian groups, with given degree and diameter. 

In order to simplify notation, we will also use $B_k(p)$ to denote the set of points with integer coordinates contained in the closed ball $B_k(p)$. The order of the largest subgraph $S$ of degree $\Delta \leq 2k$ and diameter $D = 2p$ or $D = 2p+1$, that can be constructed on the $k$-dimensional mesh, will be denoted $N_k(\Delta, p)$. Alternatively we could use the notations $B_k(D)$ and $N_k(\Delta, D)$, specifying whether $D$ is even or odd. If $k' > k$, the following inequalities are straightforward:

\begin{equation}
\label{eq:inequalities}
| B_k(p) | \leq N_{k'}(2k, p) \leq | B_{k'}(p) |
\end{equation}

\begin{figure}[htbp]
\begin{center}
	 	\includegraphics[width=0.3\textwidth]{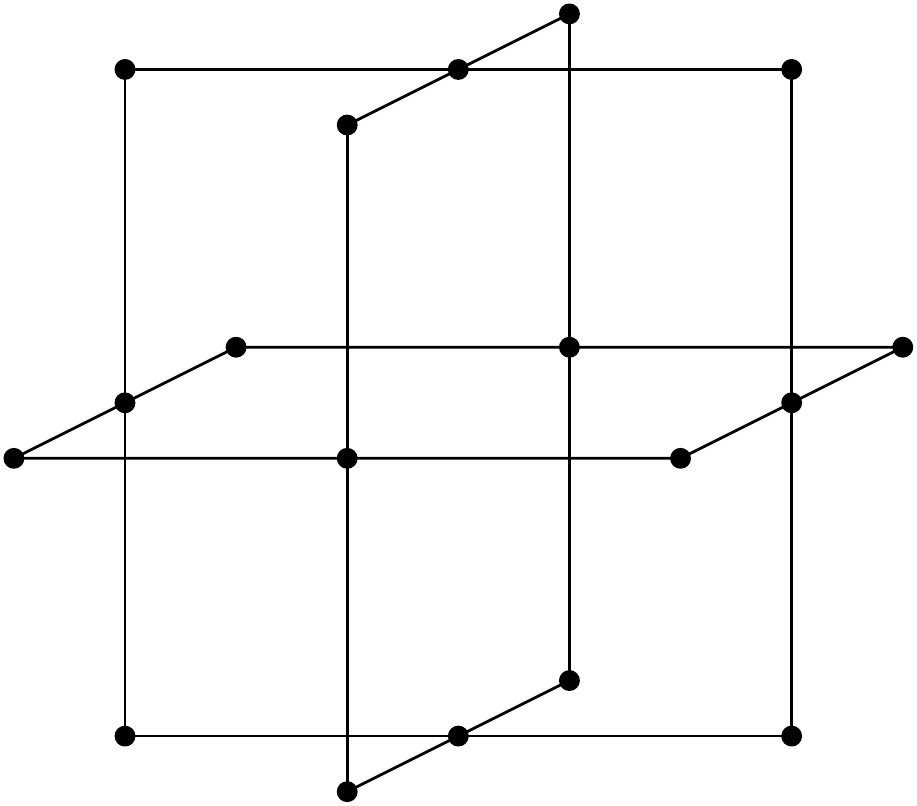}
	 \caption{Construction for $\Delta = 4$ and $D = 4$ in the 3D mesh}
	 \label{fig:example3D}
\end{center}	
\end{figure}

The first inequality tells us that if we go to a higher dimension, keeping $\Delta$ and $p$ constant, we can construct larger subgraphs.  The reason for that is that we can move along the extra dimensions in order to avoid \lq collisions\rq. Figure \ref{fig:example3D} is an example of one such construction in dimension $k = 3$, of a subgraph with degree $\Delta = 4$, diameter $D = 4$ (i.e. $p = 2$), and $18$ vertices, whereas $| B_2(2) | = 13$.

From (\ref{eq:inequalities}) follows the importance of determining the numbers $| B_{k}(p) |$. Counting the number of lattice points contained in a circle is a problem that goes back to Gauss, and there are several approximate results for the number of lattice points in balls and other sets, e.g. \cite{skriganov, widmer}. Regarding the \emph{exact} number of lattice points contained in closed balls in the $L^1$ metric in arbitrary dimension, the main reference seems to be a paper by Vassilev-Missana and Atanassov \cite{atanassov}. The following result was given in \cite{maxddbs}, and we reproduce it here with slight modifications.

\begin{theorem}
\label{th:manhattanballs}
The cardinality of $B_k(p)$ is

\begin{equation}
\label{eq:result1}
| B_k(p) | = \left\{ \begin{array}{ll}
\sum_{i=0}^{p}{{k \choose i}{k+p-i \choose p-i}} = \sum_{i=0}^{p}{{k \choose p-i}{k+i \choose i}} & \textrm{if $D = 2p$}\\
2\sum_{i=0}^{p}{{k-1 \choose i}{k+p-i \choose p-i}} = 2\sum_{i=0}^{p}{{k-1 \choose p-i}{k+i \choose i}} & \textrm{if $D = 2p+1$}
\end{array} \right.
\end{equation}
\end{theorem}

\textbf{Proof:}

\begin{figure}[htbp]
\begin{center}
	 	\includegraphics[width=0.6\textwidth]{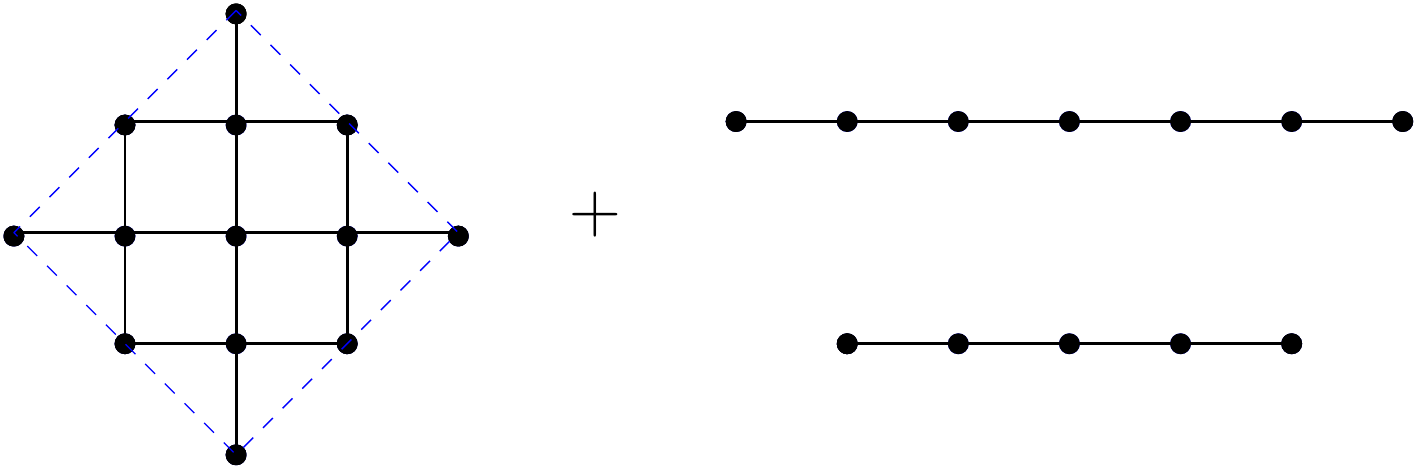}
	 \caption{Decomposition of the ball of diameter 6}
	 \label{fig:decomp}
\end{center}	
\end{figure}

The ball $B_k(p)$ can be constructed as the union of a smaller ball in the same dimension $k$, plus two balls in dimension $k-1$. Let $D$ be even, and let us place the origin of our coordinate system in the central lattice point. Then the subset consisting of all the lattice points having the $k$-th coordinate $x_k$ equal to zero is $B_{k-1}(p)$. This subset separates $B_k(p)$ into two hemispheres, one made up by those lattice points with a positive $k$-th coordinate, and those with a negative $k$-th coordinate. The layers with $x_k = \pm 1$ are $B_{k-1}(p-1)$, i.e. they have diameter $D-2 = 2(p-1)$. If we remove one of these layers (say, the one with $x_k = -1$), and put both hemispheres together, we get $B_k(p-1)$. Figure \ref{fig:decomp} shows the decomposition for $k=2$ and $p=3$. The case of odd diameter is very similar.  

From this decomposition we get the recurrence relation

\begin{equation}
\label{eq:recurrence}
f(k,p) = f(k,p-1) + f(k-1,p) + f(k-1,p-1)
\end{equation}

where $f(k,p)$ denotes the number of lattice points in $B_k(p)$. The boundary conditions are:

\begin{displaymath}
f(k,0) = \left\{ \begin{array}{ll}
1 & \textrm{if $D$ is even}\\
2 & \textrm{if $D$ is odd}
\end{array} \right.
\end{displaymath}

\begin{equation}
\label{eq:boundary2}
f(1,p) = \left\{ \begin{array}{ll}
2p+1 & \textrm{if $D$ is even}\\
2(p+1) & \textrm{if $D$ is odd}
\end{array} \right.
\end{equation}

We want to find the generating function $A_k(z) = \sum_{p \geq 0}{f(k,p)z^p}$. Multiplying (\ref{eq:recurrence}) by $z^p$ and summing over $p \geq 1$ we get

\begin{equation}
\label{eq:gen1}
A_k(z) - A_k(0) = zA_k(z) + (A_{k-1}(z)-A_{k-1}(0)) + zA_{k-1}(z)
\end{equation}

whence

\begin{equation}
\label{eq:gen2}
A_k(z) = \frac{1+z}{1-z} A_{k-1}(z)
\end{equation}

With the aid of the boundary conditions (\ref{eq:boundary2}) we get

\begin{equation}
\label{eq:gen3}
A_k(z) = \left\{ \begin{array}{ll}
\frac{(1+z)^k}{(1-z)^{k+1}} & \textrm{if $D = 2p$}\\
2\frac{(1+z)^{k-1}}{(1-z)^{k+1}} & \textrm{if $D = 2p+1$}
\end{array} \right.
\end{equation}

For even $D$, $A_k(z)$ is the product of $(1+z)^k = \sum_p{k \choose p}z^p$ and $1/(1-z)^{k+1} = \sum_p{k+p \choose p}z^p$. Then, the series of $A_k(z)$ can be obtained as the convolution of the respective factor series. The series of $A_k(z)$ for odd $D$ can be obtained in the same manner.   \EndProof

For even $D$, the numbers $| B_k(p) |$ turn out to be the \emph{Delannoy numbers} (sequence A008288 of \cite{oeis}), which appear in a variety of combinatorial and geometric problems \cite{sulanke}. This particular interpretation of Delannoy numbers was first given by Vassilev-Missana and Atanassov \cite{atanassov}, and later rediscovered by Schr\"oder \cite{schroeder}, and then by us. Our formulation and proof are different from the ones in \cite{atanassov, schroeder}. For odd $D$, the numbers $| B_k(p) |$ are known as a Riordan array of coordination sequences (sequence A113413 of \cite{oeis}). Tables \ref{tab:even} and \ref{tab:odd} show the first few values of $| B_k(p) |$ for even and odd $D$. They can be constructed in a Pascal-like fashion, with the convention that $| B_0(p) | = 1$. 

\begin{table}[htp]
\begin{center}
\begin{tabular}{|cc|*{9}{c}|} \hline
\multicolumn{2}{|c|}{} & \multicolumn{9}{c|}{$p$} \\
\multicolumn{2}{|c|}{$k$} & 0 & 1 & 2 & 3 & 4 & 5 & 6 & 7 & 8\\ \hline 
\multicolumn{2}{|c|}{0} & 1 & 1 & 1 & 1 & 1 & 1 & 1 & 1 & 1 \\ 
\multicolumn{2}{|c|}{1} & 1 & 3 & 5 & 7 & 9 & 11 & 13 & 15 & 17 \\ 
\multicolumn{2}{|c|}{2} & 1 & 5 & 13 & 25 & 41 & 61 & 85 & 113 & 145 \\ 
\multicolumn{2}{|c|}{3} & 1 & 7 & 25 & 63 & 129 & 231 & 377 & 575 & 833 \\ 
\multicolumn{2}{|c|}{4} & 1 & 9 & 41 & 129 & 321 & 681 & 1289 & 2241 & 3649 \\ \hline
\end{tabular}
\caption{Some values of $| B_k(p) |$ for even $D$}
\label{tab:even}
\end{center}
\end{table}

\begin{table}[htp]
\begin{center}
\begin{tabular}{|cc|*{9}{c}|} \hline
\multicolumn{2}{|c|}{} & \multicolumn{9}{c|}{$p$} \\
\multicolumn{2}{|c|}{$k$} & 0 & 1 & 2 & 3 & 4 & 5 & 6 & 7 & 8 \\ \hline 
\multicolumn{2}{|c|}{0} & 1 & 1 & 1 & 1 & 1 & 1 & 1 & 1 & 1 \\ 
\multicolumn{2}{|c|}{1} & 2 & 4 & 6 & 8 & 10 & 12 & 14 & 16 & 18 \\ 
\multicolumn{2}{|c|}{2} & 2 & 8 & 18 & 32 & 50 & 72 & 98 & 128 & 162 \\ 
\multicolumn{2}{|c|}{3} & 2 & 12 & 38 & 88 & 170 & 292 & 462 & 688 & 978 \\ 
\multicolumn{2}{|c|}{4} & 2 & 16 & 66 & 192 & 450 & 912 & 1666 & 2816 & 4482 \\ \hline
\end{tabular}
\caption{Some values of $| B_k(p) |$ for odd $D$}
\label{tab:odd}
\end{center}
\end{table}

It is known that Delannoy numbers have no closed form, meaning that they cannot be represented as a linear combination of a fixed number of hypergeometric terms (which can be verified with the aid of the methods developed in \cite{zeilberger}). However, we can extract asymptotic information from the generating function $A_k(z)$ in the proof of Theorem \ref{th:manhattanballs} above. Recall that $\alpha$ is an \emph{algebraic singularity} of the function $f$ if $f$ can be written near $\alpha$ as

\begin{equation}
\label{eq:singularity}
f(z) = f_0(z) + \frac{g(z)}{(1-z/\alpha)^\omega}
\end{equation}

where $f_0$ and $g$ are analytic near $\alpha$, $g$ is nonzero near $\alpha$, and $\omega$ is a real number different from $0, -1, -2, \ldots$. We readily recognize that $1$ is an algebraic singularity of $A_k(z)$, since $A_k(z)$ can be written in the above form, with $g(z) = (1+z)^k$ for $D$ even, and $g(z) = 2(1+z)^{k-1}$ for $D$ odd, and all the other conditions are satisfied. Now we can readily apply Theorem 3 of \cite{Lue80}:

\begin{theorem}
\label{th:asymptotic1}
Suppose that for some real $r > 0$, $A(z)$ is analytic in the region $| z | < r$, and has a finite number $t>0$ of singularities on the circle $| z | = r$, all of which are algebraic. Let $\alpha_i, \ \omega_i$, and $g_i$ be the values of $\alpha, \ \omega$, and $g$ in (\ref{eq:singularity}), corresponding to the $i$-th such singularity. Then $A(z)$ is the generating function for a sequence $\langle a_n \rangle$ satisfying
\begin{displaymath}
\label{eq:asymptotic1}
a_n = \frac{1}{n} \sum_{i=1}^{t} \frac{g_i(\alpha_i)n^{\omega_i}}{\Gamma(\omega_i)\alpha_i^n} + o(r^{-n}n^{\Omega-1})
\end{displaymath}
where $\Omega$ is the maximum of the $\omega_i$ and $\Gamma$ denotes the Gamma function.
\end{theorem}

We get

\begin{corollary}
\label{coro:asymptotic2}
\begin{displaymath}
\label{eq:asymptotic2}
| B_k(p) | = \frac{(2p)^k}{\Gamma(k+1)} + o(p^k) = \frac{(2p)^k}{k!} + o(p^k)
\end{displaymath}
\end{corollary}
\EndProof


\section{Subgraphs of degree $4$ in the $3$-dimensional mesh}
\label{sec:3D}

An interesting special case is $k = 3$ and $\Delta = 4$. In this case, the inequalities (\ref{eq:inequalities}) translate to $| B_2(p) | \leq N_{3}(4, p) \leq | B_{3}(p) |$, or

\begin{displaymath}
\begin{array}{ll}
2p^2 + 2p + 1 \leq N_{3}(4, p) \leq 4p^3/3 + 2p^2 + 8p/3 + 1   & \textrm{if $D = 2p$}\\
2(p^2 + 2p + 1) \leq N_{3}(4, p) \leq 4p^3/3 + 4p^2 + 14p/3 + 2   & \textrm{if $D = 2p+1$}
\end{array}
\end{displaymath}

The following theorem shows that the lower bounds are in fact a lot closer to the upper bounds:

\begin{theorem}
\label{th:newbounds3D}
\begin{equation}
\begin{array}{llll}
4p^3/3 + 2p^2 - 4p/3 + 3 & \leq \  N_{3}(4, p) & \leq \  4p^3/3 + 2p^2 + 8p/3 + 1 & \quad \textrm{if $D = 2p$}\\
4p^3/3 + 4p^2 + 2p/3 & \leq \  N_{3}(4, p) & \leq \  4p^3/3 + 4p^2 + 14p/3 + 2   & \quad \textrm{if $D = 2p+1$}
\end{array}
\end{equation}
\end{theorem}

\textbf{Proof:}

\begin{figure}[htbp]
\begin{center}
	 	\includegraphics[width=0.7\textwidth]{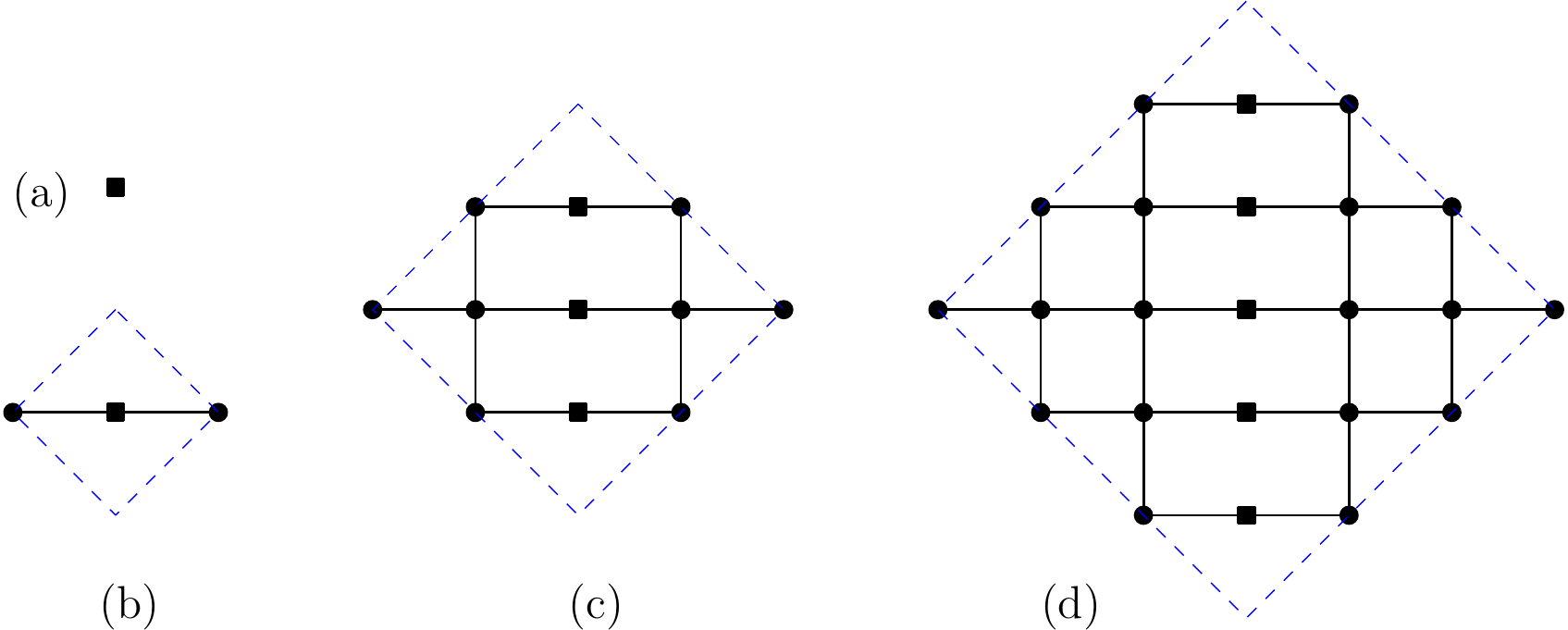}
	 \caption{Construction for $\Delta = 4$ and even diameter in the 3D mesh}
	 \label{fig:evenmesh}
\end{center}	
\end{figure}

Let $D = 2p$, with $p \geq 2$, and let us go back to our $L^1$ metric space model in dimension two. W.l.o.g. we can pick the center of our coordinate system as the center of all our balls and constructions. Now let $E_2(p)$ be the graph obtained from $B_2(p)$ by removing all the edges along the $y$-axis, and the two vertices that are left isolated (see Figure \ref{fig:evenmesh}).

Note that in $E_2(p)$ there are two classes of nodes: the ones along the $y$-axis, shown as black squares, and all the other ones, shown as black circles. The interior circular nodes have degree $4$, while the interior square nodes only have degree $2$, and we will use them to move along the third dimension $z$.

Our three-dimensional graphs $H_3(p)$ (with $p \geq 2$) will consist of layers of $E_2(i)$ connected via the square nodes: At $z=0$ we have $E_2(p)$, and at $z = \pm i$ we have $E_2(p-i)$, for $1 \leq i \leq p$. It is straightforward to check that $H_3(p)$ has diameter $D = 2p$, so let us concentrate on the number of nodes.

The number of vertices on each layer $E_2(i)$ is $| B_2(i) |$ minus two \lq missing\rq\ vertices. Therefore

\begin{displaymath}
\begin{array}{ll}
|E_2(p)| & = 2p^2 + 2p - 1 \\
|H_3(p)| & = 2 \sum_{i=1}^{p-1} (2i^2 + 2i - 1) + (2p^2 + 2p - 1) \\
         & = 4p^3/3 + 2p^2 - 4p/3 + 3
\end{array}
\end{displaymath}

\begin{figure}[htbp]
\begin{center}
	 	\includegraphics[width=0.65\textwidth]{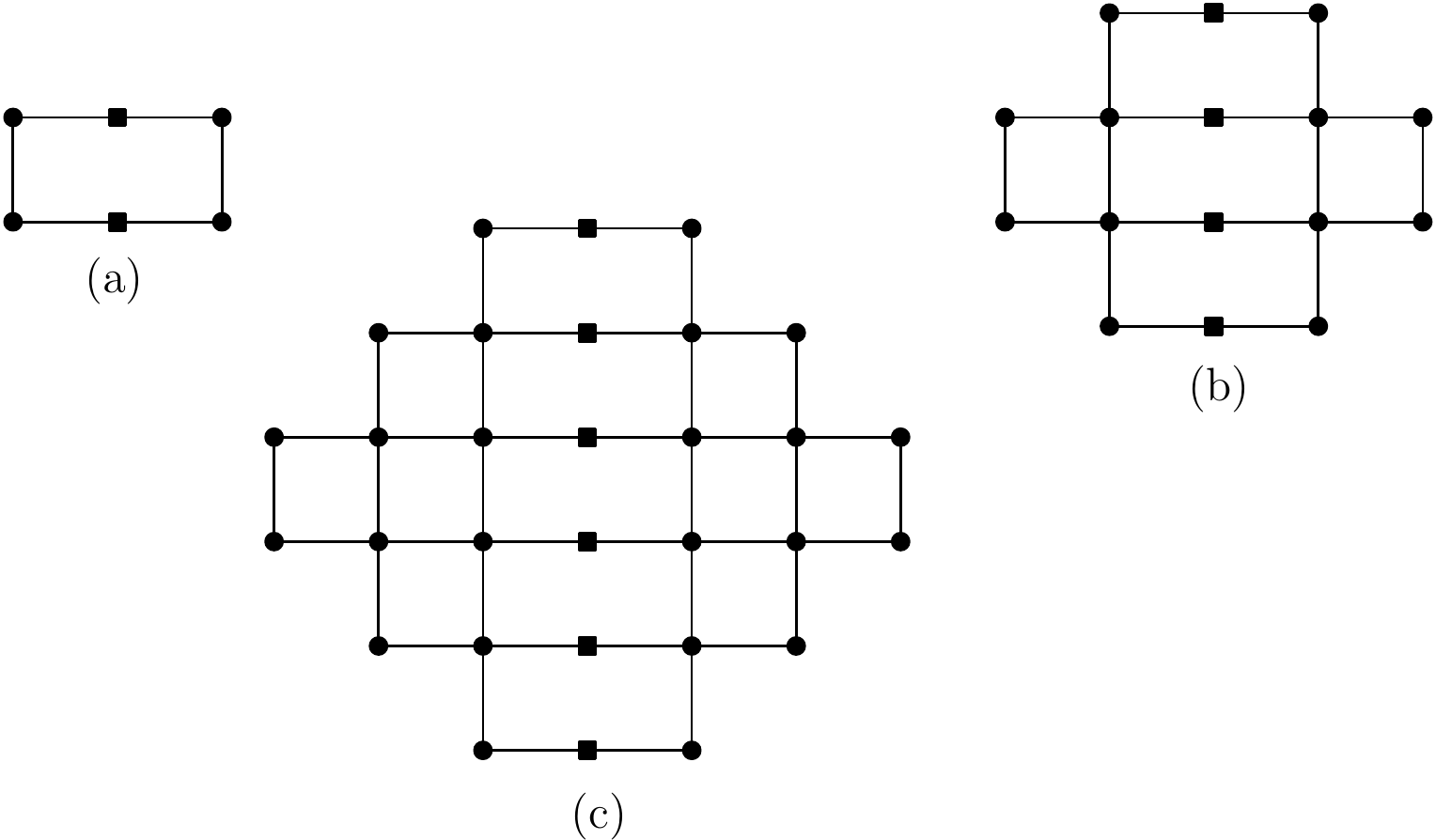}
	 \caption{Construction for $\Delta = 4$ and odd diameter in the 3D mesh}
	 \label{fig:oddmesh}
\end{center}	
\end{figure}

Now let $D = 2p+1$, with $p \geq 2$. The construction here is also made with layers that are $B_2(p)$ with the central \lq spine\rq\ suppressed (i.e. the edges on the $y$-axis and the two tip vertices). Let us call these graphs $O_2(p)$. Figure \ref{fig:oddmesh} shows $O_2(1)$, $O_2(2)$, and $O_2(3)$. The three-dimensional graph $Q_3(p)$ has a layer of type $O_2(p-i)$ at $z = \pm i$, for $0 \leq i \leq p-1$.

Again, it is easy to check that the diameter is correct. Regarding the number of nodes, in every $O_2(i)$ there are only two nodes missing from $B_2(i)$. Hence

\begin{displaymath}
\begin{array}{ll}
|Q_3(p)| & = 4 \sum_{i=1}^{p-1} (i^2 + 2i) + 2p^2+4p \\
         & = (4p^3 + 12p^2 + 2p)/3
\end{array}
\end{displaymath}
\EndProof

Note that the above constructions are asymptotically optimal, since they agree up to the second term with the upper bounds. They are also optimal in another sense, as shown by the following 

\begin{corollary}
\label{coro:edges1}
The average degree $\hat{\delta}$ of the constructions in Theorem \ref{th:newbounds3D} tends to $4$ as $p$ approaches infinity. 
\end{corollary}

\textbf{Proof:}

$\hat{\delta} = 2A/V$, where $A$ is the number of edges and $V$ is the number of vertices. The number of edges is $8p^3/3+4p/3+2$ for $D=2p$, and  $8p^3/3+4p^2-2p/3$ for $D=2p+1$. 
\EndProof

It is also very likely that these constructions can be extended to higher dimensions.


\section{Subgraphs of degree $3$ in the $2$-dimensional mesh}
\label{sec:2D}

The smallest case of $\Delta < 2k$ that makes sense is $\Delta = 3$ in dimension two. In this case we will also show that the lower bounds are quite close to the upper bounds. We have the following

\begin{theorem}
\label{th:newbounds2D}
\begin{equation}
\begin{array}{llll}
2p^2-2p+1 & \leq \  N_{2}(3, p) & \leq \  2p^2+2p+1 & \quad \textrm{if $D = 2p$}\\
8r^2+2r = 2p^2+p & \leq \  N_{2}(3, p) & \leq \  2p^2+4p+2   & \quad \textrm{if $D = 4r+1$}\\
8r^2+10r+6 = 2p^2+p+3 & \leq \  N_{2}(3, p) & \leq \  2p^2+4p+2   & \quad \textrm{if $D = 4r+3$}
\end{array}
\end{equation}
\end{theorem}

\textbf{Proof:}

As in Section \ref{sec:3D}, we will give two constructions that achieve the new bounds. In this case, a rigorous mathematical description of our constructions will be quite cumbersome; instead, we will use some geometric analogies to describe them. What we do is that we take the balls $B_{2}(p)$ and we try to fill them up in a convenient way with as many \lq building blocks\rq \ as possible. Our construction elements will be $2 \times 2$ square blocks, and $2 \times 1$ half-blocks, or rectangular \lq bricks\rq. We have to pack them up in such a way that no four edges meet at a single point.

\begin{figure}[htbp]
\begin{center}
	 	\includegraphics[width=0.8\textwidth]{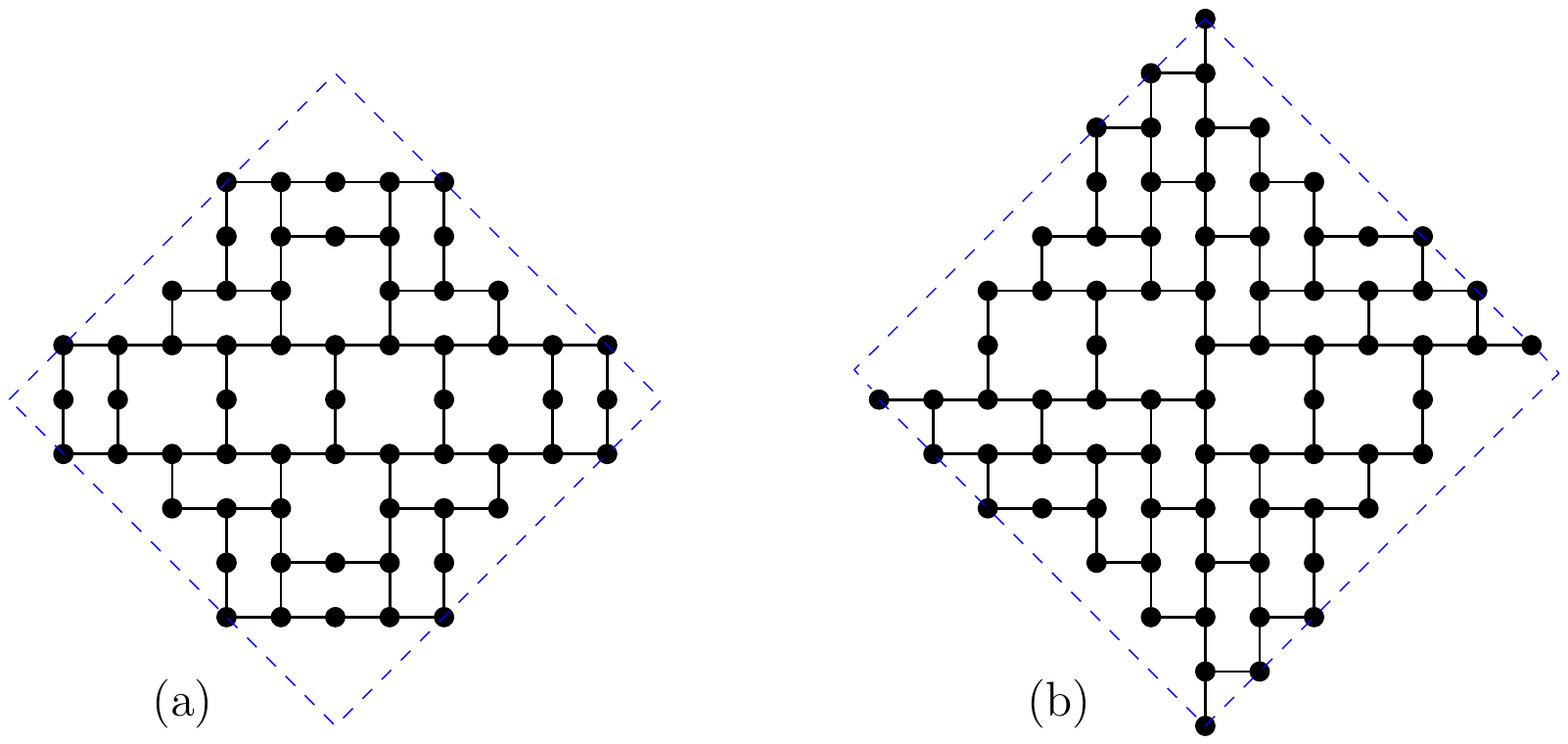}
	 \caption{Constructions for $\Delta = 3$ in the two-dimensional mesh}
	 \label{fig:delta3}
\end{center}	
\end{figure}

Let us start with even $D$: We place as many square blocks as possible along the $x$ and $y$ axes, and if there is some space left at the end, we fill it up with bricks, as shown in Figure \ref{fig:delta3} (a). On each horizontal semi-axis we will have to use $(p-1)/2$ square blocks, and on each vertical semi-axis we will need $(p-3)/2$ square blocks.\footnote{A fraction here means that we have to complete with a brick at the end.}

\begin{figure}[htbp]
\begin{center}
	 	\includegraphics[width=0.25\textwidth]{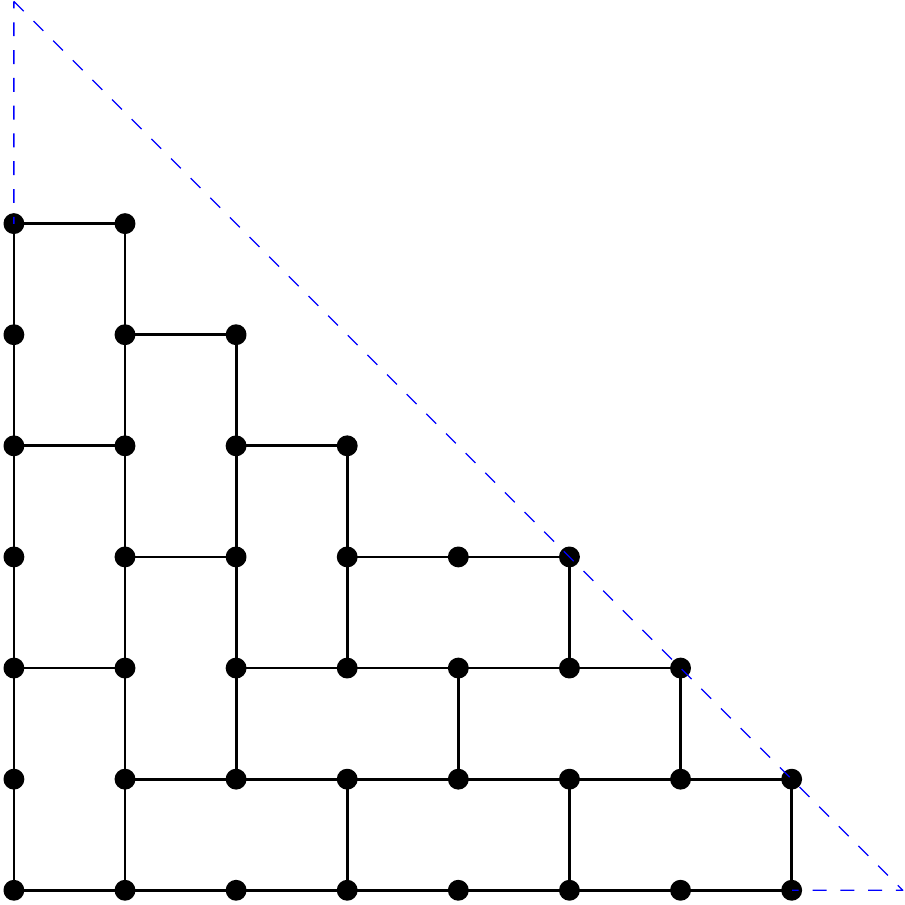}
	 \caption{Brick pattern in a triangular region}
	 \label{fig:triangle}
\end{center}	
\end{figure}

Now we will be left with four empty triangular regions that will have to be occupied. They are isosceles triangles, with side $p-2$. We will fill them with bricks, forming a pattern of interwoven horizontal and vertical bricks, as shown in Figure \ref{fig:triangle}. Starting from the innermost corner there is only one way to do this, since every other vertex on the triangle sides already has three edges attached to it. Moreover, all triangles constructed this way will be the same, save rotations and reflections. We will make use of this fact again in the case of odd $D$. The construction ensures that no vertex has degree $4$.

Next we have to show that the diameter is $D = 2p$. First consider the blocks along the axes. We shall call the vertices belonging to one of those blocks \lq block vertices\rq. It is clear that starting from any block vertex we can reach any other block vertex in at most $2p$ steps. Hence we only have to include the vertices in the triangular regions, and show that for any such vertex, we can reach any other vertex in at most $2p$ steps. The procedure is quite simple: First note that we can reach the innermost vertex of the triangle (which is also a block vertex) in at most $p-2$ steps (the length of the triangle sides). Call this vertex a \lq corner\rq. From there we can reach any block vertex in at most $p+2$ steps, and we can reach any other corner in at most $4$ steps, which means that we can reach any other vertex located in a triangular region in at most $2p$ steps, as desired.

Now we are going to count the number of vertices: If $p$ is even we have $2(p-3)$ square blocks, plus four bricks, one at the end of each semi-axis. The square blocks contribute $10(p-3)+3$ points, and we have to add $12$ additional points corresponding to the four bricks. Each triangle contributes $2\lfloor (p-2)/2 \rfloor^2 - 2\lfloor (p-2)/2 \rfloor = [(p-3)^2-1]/2$ points. Hence, the total number of vertices is $2p^2-2p+1$. If $p$ is odd we have $2(p-2)$ square blocks, and no bricks at the end of the semi-axes. The square blocks contribute $10(p-2)+3$ points. Finally, each triangle contributes $2\lfloor (p-2)/2 \rfloor^2 = (p-3)^2/2$ points. Hence, the total number of vertices is also $2p^2-2p+1$.

Now consider the construction for odd $D$. In this case we have a central \lq spine\rq \ of length $D = 2p+1$ (the $y$-axis), and two horizontal semi-axes of length $p$. We place $\lceil (p-2)/2 \rceil$ square blocks on each horizontal semi-axis, as shown in Figure \ref{fig:delta3} (b). Again, we are left with four triangular regions that are filled up with interwoven bricks, in the same manner as before. In this case, two triangles have sides of length $p$, and the other two have sides of length $p-1$.

We show that the diameter is $D = 2p+1$. We proceed as in the case of even $D$. In this case we distinguish block vertices and \lq axis vertices\rq \ (i.e. vertices belonging to the central spine or to the horizontal semi-axes). We also distinguish two central vertices at the intersection of the vertical spine and the two horizontal semi-axes. It is clear that starting from any block or axis vertex we can reach the closest central vertex in at most $p$ steps. Now, if we start at a vertex belonging to one of the larger triangles, of side length $p$, we can reach the corresponding corner (which is one of the central vertices) in at most $p$ steps. Starting at a vertex in one of the smaller triangles, of side length $p-1$, we can reach the corresponding corner in at most $p-1$ steps, and the nearest central vertex is just one step away. Since both central vertices are adjacent, the proof of the diameter follows.

Regarding the number of vertices, we have $2p+2$ vertices on the vertical spine, plus $2p$ vertices from the horizontal semi-axes, plus $6\lceil (p-2)/2 \rceil$ block vertices, plus the triangles. If $p$ is even, the two large triangles contribute $p(p-2)$ extra vertices, and the two smaller triangles contribute $(p-2)^2$ extra vertices. If $p$ is odd, the two large triangles contribute $(p-1)^2$ extra vertices, and the two smaller triangles contribute $(p-1)(p-3)$ extra vertices. Hence the total number of vertices is $2p^2+p$ if $p$ is even, and $2p^2+p+3$ if $p$ is odd.

\EndProof

These constructions are also asymptotically optimal in the sense of the average degree: 

\begin{corollary}
\label{coro:edges2}
The average degree $\hat{\delta}$ of the constructions in Theorem \ref{th:newbounds2D} tends to $3$ as $p$ approaches infinity. 
\end{corollary}

\textbf{Proof:}

The number of edges is 

\begin{equation}
\begin{array}{llll}
3p^2-6p+6 & \quad \textrm{if $D = 4r$} \\
3p^2-p-1 & \quad \textrm{if $D = 4r+1$}\\
3p^2-6p+4 & \quad \textrm{if $D = 4r+2$} \\
3p^2-p+3 & \quad \textrm{if $D = 4r+3$}\\
\end{array}
\end{equation}
\EndProof

Despite being asymptotically optimal, these constructions can be improved without much effort for particular values of $D$, but extracting a general pattern and counting the vertices for arbitrary $D$ may prove a challenge. In Table \ref{tab:orders} we give the order of the largest graphs that have been constructed for some small values of $D$. The actual graphs can be seen in \cite{combinatoricswiki}. Note that for diameters $D = 2, \ldots, 6$, the upper bound given in the table is smaller than the upper bound given by Theorem \ref{th:newbounds2D}, and the largest graphs are optimal. These graphs have been obtained by exhaustive computer search. 

\begin{table}[htp]
\begin{center}
\begin{tabular}{|cc|*{15}{c}|} \hline
\multicolumn{2}{|c|}{Diameter} & 2 & 3 & 4 & 5 & 6 & 7 & 8 & 9 & 10 & 11 & 12 & 13 & 14 & 15 & 16\\ \hline
\multicolumn{2}{|c|}{Upper bound} & 4 & 6 & 10 & 14 & 22 & 32 & 41 & 50 & 61 & 72 & 85 & 98 & 113 & 128 & 145\\ 
\multicolumn{2}{|c|}{Order of largest graph} & 4 & 6 & 10 & 14 & 22 & 28 & 37 & 44 & 52 & 68 & 77 & 90 & 104 & 124 & 135\\ \hline
\end{tabular}
\caption{Orders of the largest known subgraphs of the two-dimensional mesh, with $\Delta = 3$ and $2 \leq D \leq 16$.}
\label{tab:orders}
\end{center}
\end{table}


\section{Conclusions and open problems}
\label{sec:open}

To us, it is always amazing that a simple combinatorial setting like this one provides such a wealth of interesting and difficult problems. In the preceding two sections we have seen constructions that result in lower bounds for $N_{k}(\Delta, p)$, that are  asymptotically close to the upper bounds, for small values of $k$ and $\Delta$, with $\Delta < 2k$. Now, is that also possible for arbitrary $k$ and $\Delta < 2k$? How close can we get to the upper bounds? The constructions in Section \ref{sec:3D} could be generalized to higher dimensions, more precisely to the case $\Delta = 2k-2$. The constructions in Section \ref{sec:2D} might also be extended to $\Delta = 5$ in dimension $3$, but it does not appear to be easy.

There are several additional properties that could be considered for every construction, like connectivity ({\it sic} fault-tolerance), average path length, symmetry, etc. For example, in parallel computing applications, our constructions of Section \ref{sec:3D} might not be the best ones, as they impose a relatively high communication overhead on the square vertices.

On the other hand, what is the computational complexity of finding the largest degree\&diameter bounded subgraph in the mesh? Does it remain $\NP-$hard for all dimensions?, for some dimensions?, for any dimension? If it remains $\NP-$hard, then, can it be approximated to within a constant ratio?

This same study can be carried out for other host networks of theoretical and/or practical importance: the hypercube, the butterfly, the cube-connected cycles, Cayley graphs, etc. In \cite{maxddbs} there is a preliminary discussion for hypercubes, and a very incomplete heuristic study was attempted for some random networks, but other than that, the field remains totally virgin. 

Finally, it is worth noting that this problem could be related to a well-known problem arising in parallel computing, namely that \lq embedding\rq \ or \lq emulating\rq \ an arbitrary graph in a mesh (see \cite{zienicke}, for instance). This connection was indicated to us by one of the reviewers, and it is also interesting to explore further.

\section*{Acknowledgements}

We wish to thank Charles Delorme, Robert Israel, and Robert Sulanke, for their helpful comments and suggestions, and Doron Zeilberger and Marko Petkov\v{s}ek for their help with the algorithm \lq Hyper\rq, to verify that Delannoy numbers have no closed form. Krassimir Atanassov kindly agreed to share his manuscript about Delannoy numbers and lattice points, and provided other important references. Anthony Dekker obtained the optimal graphs for diameters $2, \ldots, 6$ in Table \ref{tab:orders} by exhaustive computer search. The second author wishes to thank his nine-year old daughter Karla, who also helped obtain some of the graphs listed in Table \ref{tab:orders}. Finally, we are indebted to the anonymous referees who reviewed the paper, for their constructive and useful suggestions, which not only helped us improve the paper, but also contain guidelines for our future research.


\end{document}